\documentclass[12pt]{article}
\usepackage{mathrsfs}
\usepackage{amssymb}
\usepackage{amsmath,amsthm}

\newtheorem{Theorem}{Theorem}
\newtheorem{Lemma}[Theorem]{Lemma}
\newtheorem{Corollary}[Theorem]{Corollary}

\theoremstyle{definition}
\newtheorem*{pf}{Proof}

\begin{document}

\title{A mean value formula for elliptic curves}
\author{Rongquan Feng$^{1}$, Hongfeng Wu$^{2}$
\\
{\small $^1$LMAM, School of Mathematical Sciences, Peking University, Beijing 100871,  China}\\
{\small\small fengrq@math.pku.edu.cn}\\[1ex]
{\small $^2$College of Science, North China University of technology, Beijing 100144, P.R. China}\\
{\small\small whfmath@gmail.com}}

\date{}
\maketitle

\begin{abstract}
It is proved in this paper that for any point on an elliptic curve,
the mean value of $x$-coordinates of its $n$-division points is the
same as its $x$-coordinate and that of $y$-coordinates of its
$n$-division points is $n$ times of its $y$-coordinate.
\end{abstract}


\vspace*{1cm}

\section{Introduction}
Let $K$ be a field with $\text{char}(K)\neq 2, 3$ and let
$\overline{K}$ be the algebraic closure of $K$. Every elliptic curve
$E$ over $K$ can be written as a classical Weierstrass equation
$$E: y^2=x^3+ax+b$$ with coefficients
$a,\,b\in K$. A point $Q$ on $E$ is said to be smooth (or
non-singular) if $\left(\frac{\partial f}{\partial
x}|_Q,\frac{\partial f}{\partial y}|_Q\right)\ne (0,0)$, where
$f(x,y)=y^2-x^3-ax-b$. The point multiplication is the operation of
computing $$nP=\underbrace{P+P+\cdots+P}_{n}$$ for any point $P\in
E$ and a positive integer $n$. The multiplication-by-$n$ map
$$\begin{array}{rcll}[n]:&E&\rightarrow &E\\
&P&\mapsto &nP\end{array}$$ is an isogeny of degree $n^2$. For a
point $Q\in E$, any element of $[n]^{-1}(Q)$ is called an
$n$-division point of $Q$. Assume that $(\text{char}(K),n)=1$. In
this paper, the following result on the mean value of the $x,
y$-coordinates of all the $n$-division points of any smooth point on
an elliptic curve is proved.

\begin{Theorem}\label{meanvalue}
Let $E$ be an elliptic curve defined over $K$, and let
$Q=(x_Q,y_Q)\in E$ be a point with $Q\neq \mathcal{O}$. Set
$$\Lambda=\{P=(x_P,y_P)\in E(\overline{K})\mid nP=Q\}.$$ Then $$\frac{1}{n^2}\sum\limits_{P\in
\Lambda}x_P=x_Q$$
and
$$\frac{1}{n^2}\sum\limits_{P\in\Lambda}y_P=ny_Q.$$
\end{Theorem}
~\\

According to Theorem \ref{meanvalue}, let $P_i=(x_i,y_i),
i=1,2,\cdots,n^2$, be all the points such that $nP=Q$ and let
$\lambda_i$ be the slope of the line through $P_i$ and $Q$, then
$y_Q=\lambda_i(x_Q-x_i)+y_i$. Therefore,
$$n^2y_Q=\sum\limits_{i=1}^{n^2}\lambda_i\cdot
(\sum\limits_{i=1}^{n^2}x_i)/n^2-\sum\limits_{i=1}^{n^2}\lambda_ix_i+\sum\limits_{i=1}^{n^2}y_i.$$
Thus we have
$$y_Q=\dfrac{\sum\limits_{i=1}^{n^2}\lambda_i}{n^2}\cdot
\dfrac{\sum\limits_{i=1}^{n^2}x_i}{n^2}-\dfrac{\sum\limits_{i=1}^{n^2}\lambda_ix_i}{n^2}+\dfrac{\sum\limits_{i=1}^{n^2}y_i}{n^2}
=\overline{\lambda_i}\cdot
\overline{x_i}-\overline{\lambda_ix_i}+\overline{y_i},$$ where
$\overline{\lambda_i},~
\overline{x_i},~\overline{\lambda_ix_i},~\overline{y_i}$ are the
average values of the variables $\lambda_i, x_i, \lambda_ix_i$ and
$y_i$, respectively. Therefore,
$$Q=(x_Q,y_Q)=(\overline{x_i},~\overline{\lambda_i}\cdot
\overline{x_i}-\overline{\lambda_ix_i}+\overline{y_i})=\left(\overline{x_i},~\frac{1}{n}\overline{y_i}\right).$$

\noindent{\bf Remark:} The discrete logarithm problem in elliptic
curve $E$ is to find $n$ by given $P,\,Q\in E$ with $Q=nP$. The
above theorem gives some information on the integer $n$.

\section{Proof of Theorem \ref{meanvalue}}
To prove Theorem \ref{meanvalue}, define division polynomials
\cite{Silverman} $\psi_n\in \mathbb{Z}[x,y,a,b]$ on an elliptic
curve $E: y^2=x^3+ax+b$, inductively as follows:
\begin{equation*}
\begin{array}{rcl}
\psi_0 &=& 0,\\[.5ex]
\psi_1 &=& 1,\\[.5ex]
\psi_2 &=& 2y,\\[.5ex]
\psi_3 &=& 3x^4+6ax^2+12bx-a^2,\\[.5ex]
\psi_4 &=& 4y(x^6+5ax^4+20bx^3-5a^2x^2-4abx-8b^2-a^3),\\[.5ex]
\psi_{2n+1} &=& \psi_{n+2}\psi_{n}^3-\psi_{n-1}\psi_{n+1}^3,
\text{~for~} n\geq 2,\\[.5ex]
2y\psi_{2n} &=&
\psi_n(\psi_{n+2}\psi_{n-1}^2-\psi_{n-2}\psi_{n+1}^2), \text{~for~}
n\geq 3.
\end{array}
\end{equation*}
 It can be checked easily by induction that the $\psi_{2n}$'s are polynomials. Moreover, $\psi_n\in
 \mathbb{Z}[x,y^2,a,b]$ when $n$ is odd, and $(2y)^{-1}\psi_n\in
 \mathbb{Z}[x,y^2,a,b]$ when $n$ is even.
 Define the polynomial $$\phi_n=x\psi_n^2-\psi_{n-1}\psi_{n+1}$$ for $n\ge 1$. Then
 $\phi_n\in \mathbb{Z}[x,y^2,a,b]$. Since $y^2=x^3+ax+b$, replacing $y^2$ by $x^3+ax+b$, one have that $\phi_n\in
 \mathbb{Z}[x,a,b]$. So we can denote it by $\phi_n(x)$. Note that,
 $\psi_n\psi_m\in \mathbb{Z}[x,a,b]$ if $n$ and $m$ have the same
 parity. Furthermore, the division polynomials $\psi_n$ have the
 following properties.

\begin{Lemma}\label{firsttwoxiang}
$$\psi_n=nx^{\frac{n^2-1}{2}}+\frac{n(n^2-1)(n^2+6)}{60}ax^{\frac{n^2-5}{2}}+\text{~lower degree terms},$$ when $n$ is odd,
and
$$\psi_n=ny\left(x^{\frac{n^2-4}{2}}+\frac{(n^2-1)(n^2+6)-30}{60}ax^{\frac{n^2-8}{2}}+\text{~lower degree terms}\right),$$
when $n$ is even.
\end{Lemma}
\begin{pf}
We prove the result by induction on $n$. It is true for $n<5$.
Assume that it holds for all $\psi_m$ with $m<n$. We give the proof
only for the case for odd $n\ge 5$. The case for even $n$ can be
proved similarly. Now let $n=2k+1$ be odd, where $k\ge 2$. If $k$ is
even, then by induction,
\begin{equation*}
\begin{array}{rcl}
\psi_k &=& ky(x^{\frac{k^2-4}{2}}+\frac{(k^2-1)(k^2+6)-30}{60}ax^{\frac{k^2-8}{2}}+\cdots),\\~\\
\psi_{k+2} &=& (k+2)y(x^{\frac{k^2+4k}{2}}+\frac{(k^2+4k+3)(k^2+4k+10)-30}{60}ax^{\frac{k^2+4k-4}{2}}+\cdots),\\~\\
\psi_{k-1} &=& (k-1)x^{\frac{k^2-2k}{2}}+\frac{(k-1)(k^2-2k)(k^2-2k+7)}{60}ax^{\frac{k^2-2k-4}{2}}+\cdots,\\~\\
\psi_{k+1} &=&
(k+1)x^{\frac{k^2+2k}{2}}+\frac{(k+1)(k^2+2k)(k^2+2k+7)}{60}ax^{\frac{k^2+2k-4}{2}}+\cdots,
\end{array}
\end{equation*}
By substituting $y^4$ by $(x^3+ax+b)^2$, we have
$$\psi_{k+2}\psi_{k}^3=k^3(k+2)\left(x^{2k^2+2k}+\frac{4(k+1)(k^3+k^2+10k+3)}{60}ax^{2k^2+2k-2}+\cdots\right),$$
and
$$\psi_{k-1}\psi_{k+1}^3=(k-1)(k+1)^3x^{2k^2+2k}+\frac{4k(k-1)(k^3+2k^2+11k+7)(k+1)^3}{60}ax^{2k^2+2k-2}+\cdots.$$
Therefore
\begin{equation*}
\begin{array}{rcl}
\psi_{2k+1}&=&\psi_{k+2}\psi_{k}^3-\psi_{k-1}\psi_{k+1}^3\\~\\
&=& (2k+1)x^{2k^2+2k}+\frac{(2k+1)(4k^2+4k)(4k^2+4k+7)}{60}ax^{2k^2+2k-2}+\cdots\\~\\
&=&
(2k+1)x^{\frac{(2k+1)^2-1}{2}}+\frac{(2k+1)((2k+1)^2-1)((2k+1)^2+6)}{60}ax^{\frac{(2k+1)^2-5}{2}}+\cdots
\end{array}
\end{equation*}
The case when $k$ is odd can be proved similarly. \qed
\end{pf}

The following corollary follows immediately from Lemma
\ref{firsttwoxiang}.

\begin{Corollary}\label{firsttwoxiangpsi}
$$\psi_n^2=n^2x^{n^2-1}-\frac{n^2(n^2-1)(n^2+6)}{30}ax^{n^2-3}+\cdots,$$
and
$$\phi_n=x^{n^2}-\frac{n^2(n^2-1)}{6}ax^{n^2-2}+\cdots.$$\qed
\end{Corollary}

\noindent {\bf Proof of Theorem \ref{meanvalue}:} Define $\omega_n$
as
$$4y\omega_n=\psi_{n+2}\psi_{n-1}^2-\psi_{n-2}\psi_{n+1}^2.$$ Then
for any $P=(x_P,y_P)\in E$, we have (\cite{Silverman})
$$nP=\left(\dfrac{\phi_n(x_P)}{\psi_n^2(x_P)},\dfrac{\omega_n(x_P,y_P)}{\psi_n(x_P,y_P)^3}\right).$$
If $nP=Q$, then $\phi_n(x_P)-x_Q\psi_n^2(x_P)=0$. Therefore, for any
$P\in \Lambda$, the $x$-coordinate of $P$ satisfies the equation
$\phi_n(x)-x_Q\psi_n^2(x)=0$. From Corollary \ref{firsttwoxiangpsi},
we have that $$\phi_n(x)-x_Q\psi_n^2(x)=x^{n^2}-n^2x_Qx^{n^2-1}+
\mbox{lower degree terms}.$$ Since $\sharp \Lambda=n^2$, every root
of $\phi_n(x)-x_Q\psi_n^2(x)$ is the $x$-coordinate of some
$P\in\Lambda$. Therefore $$\sum\limits_{P\in \Lambda}x_P=n^2x_Q$$ by
Vitae's Theorem.

Now we prove the mean value formula for $y$-coordinates. Let $K$ be
the complex number field $\mathbb{C}$ first and let $\omega_1$ and
$\omega_2$ be complex numbers which are linearly independent over
$\mathbb{R}$. Define the lattice $$L = \mathbb{Z}\omega_1 +
\mathbb{Z}\omega_2 = \{n_1\omega_1 + n_2\omega_2 \mid n_1, n_2\in
\mathbb{Z}\},$$ and the Weierstrass $\wp$-function by
$$\wp(z)=\wp(z,L)=\frac{1}{z}+\sum\limits_{\omega\in L,\omega\neq 0}\left(\frac{1}{(z-\omega)^2}-\frac{1}{\omega^2}\right).$$
For integers $k\geq 3$, define the Eisenstein series $G_k$
by
$$G_k=G_k(L)=\sum\limits_{\omega\in L,\omega\neq 0}\omega^{-k}.$$
Set $g_2=60G_4$ and $g_3=140G_6$, then
$$\wp^{'}(z)^2=4\wp(z)^3-g_2\wp(z)-g_3.$$ Let
$E$ be the elliptic curve given by $y^2=4x^3-g_2x-g_3$. Then the
map
$$
\begin{array}{rcl}
\mathbb{C}/L &\rightarrow& E(\mathbb{C}) \\
z &\mapsto& \left(\wp(z),\wp^{'}(z)\right),\\
0 &\mapsto& \infty,
\end{array}
$$
is an isomorphism of groups $\mathbb{C}/L$ and $E(\mathbb{C})$.
Conversely, it is well known \cite{Silverman} that for any elliptic
curve $E$ over $\mathbb{C}$ defined by $y^2=x^3+ax+b$, there is a
lattice $L$ such that $g_2(L)=-4a, g_3(L)=-4b$ and there is an
isomorphism between groups $\mathbb{C}/L$ and $E(\mathbb{C})$ given
by $z \mapsto \left(\wp(z),\frac{1}{2}\wp^{'}(z)\right)$ and $0 \mapsto
\infty$. Therefore, for any point $(x,y)\in E(\mathbb{C})$, we have
$(x,y)=\left(\wp(z), \frac{1}{2}\wp^{'}(z)\right)$ and $n(x,y)=\left(\wp(nz),
\frac{1}{2}\wp^{'}(nz)\right)$ for some $z\in \mathbb{C}$.

Let $Q=\left(\wp(z_Q), \frac{1}{2}\wp^{'}(z_Q)\right)$ for a
$z_Q\in\mathbb{C}$. Then for any $P_i\in \Lambda$, $1\leq i\leq
n^2$, there exist integers $j,k$ with $0\leq j,k\leq n-1$, such that
$$P_i=\left(\wp\left(\frac{z_Q}{n}+\frac{j}{n}\omega_1+\frac{k}{n}\omega_2\right),
\frac{1}{2}\wp^{'}\left(\frac{z_Q}{n}+\frac{j}{n}\omega_1+\frac{k}{n}\omega_2\right)\right).$$
Thus
$$\sum\limits_{j,k=0}^{n-1}\wp\left(\frac{z_Q}{n}+\frac{j}{n}\omega_1+\frac{k}{n}\omega_2\right)=n^2\wp(z_Q)$$
which comes from  $\sum\limits_{i=1}^{n^2}x_i=n^2x_Q$. Differential
for $z_Q$, we have
$$\sum\limits_{j,k=0}^{n-1}\wp^{'}\left(\frac{z_Q}{n}+\frac{j}{n}\omega_1+\frac{k}{n}\omega_2\right)=n^3\wp{'}(z_Q).$$
That is $$\sum\limits_{i=1}^{n^2}y_i=n^3y_Q.$$

Secondly, let $K$ be a field of characteristic $0$ and let $E$ be
the elliptic curve over $K$ given by the equation $y^2 = x^3 + ax +
b$. Then all of the equations describing the group law are defined
over $\mathbb{Q}(a,b)$. Since $\mathbb{C}$ is algebraically closed
and has infinite transcendence degree over $\mathbb{Q}$,
$\mathbb{Q}(a,b)$ can be considered as a subfield of $\mathbb{C}$.
Therefore  we can regard $E$ as an elliptic curve defined over
$\mathbb{C}$. Thus the result follows.

At last assume that $K$ is a field of characteristic $p$. Then the
elliptic curve can be viewed as one defined over some finite field
$\mathbb{F}_q$, where $q=p^m$ for some integer $m$. Without loss of
generality, let $K=\mathbb{F}_q$ for convenience. Let
$K'=\mathbb{Q}_q$ be an unramified extension of the $p$-adic numbers
$\mathbb{Q}_p$ of degree $m$, and let $\overline{E}$ be an elliptic
curve over $K'$ which is a lift of $E$. Since $(n, p)=1$, the
natural reduction map $\overline{E}[n]\rightarrow E[n]$ is an
isomorphism. Now for any point $Q\in E$ with $Q\neq \mathcal{O}$, we
have a point $\overline{Q}\in \overline{E}$ such that the reduction
point is $Q$. For any point $P_i\in E(\overline{K})$ with $nP_i=Q$,
its lifted point $\overline{P}_i$ satisfies
$n\overline{P}_i=\overline{Q}$ and $\overline{P}_i\neq
\overline{P}_j$ whenever $P_i\neq P_j$. Thus
$$\sum\limits_{i=1}^{n^2}y(\overline{P}_i)=n^3y(\overline{Q})$$ since $K'$ is a
field of characteristic 0. Therefore the formula
$\sum\limits_{i=1}^{n^2}y_i=n^3y_Q$ holds by the reduction from
$\overline{E}$ to $E$.\qed\\

\noindent{\bf Remark:} \begin{enumerate}
\item[(1)] The result for $x$-coordinate of Theorem \ref{meanvalue} holds also for the elliptic curve defined by the
general Weierstrass equation $y^2+a_1xy+a_3y=x^3+a_2x^2+a_4x+a_6$.
\item[(2)] The mean value formula for $x$-coordinates was given in the first version of
this paper \cite{Feng} with a slightly complicated proof.  The
formula for $y$-coordinates was conjectured by D. Moody based on
\cite{Feng} and numerical examples in a personal email communication
\cite{Moody}.
\item[(3)] Recently, some mean value formulae for twisted
Edwards curves \cite{Bernstein,Edwards} and other alternate models of elliptic curves were given by \cite{Moody2} and \cite{Moody3}.\end{enumerate}

\section{An application}

Let $E$ be an elliptic curve over $K$ given by the Weierstrass
equation $y^2=x^3+ax+b.$ Then we have a non-zero invariant
differential $\omega=\frac{dx}{y}$. Let $\phi\in \text{End}(E)$ be a
nonzero endomorphism. Then
$\phi^{*}\omega=\omega\circ\phi=c_{\phi}\omega$ for some
$c_{\phi}\in\overline{K}(E)$ since the space $\Omega_E$ of
differential forms on $E$ is a $1$-dimensional
$\overline{K}(E)$-vector space. Since $c_{\phi}\neq 0$ and
$\text{div}(\omega)=0$, we have
$$\text{div}(c_{\phi})=\text{div}(\phi^{*}\omega)-\text{div}(\omega)=\phi^{*}\text{div}(\omega)-\text{div}(\omega)=0.$$
Hence $c_{\phi}$ has neither zeros nor poles and
$c_{\phi}\in\overline{K}$. Let $\varphi$ and $\psi$ be two nonzero
endomorphisms, then
$$c_{\varphi+\psi}\omega=(\varphi+\psi)^{*}\omega=\varphi^{*}\omega+\psi^{*}\omega=c_{\varphi}\omega+c_{\psi}\omega
=(c_{\varphi}+c_{\psi})\omega.$$ Therefore,
$c_{\varphi+\psi}=c_{\varphi}+c_{\psi}$. For any nonzero
endomorphism $\phi$, set $\phi(x,y)=(R_{\phi}(x), yS_{\phi}(x))$,
where $R_\phi$ and $S_\phi$ are rational functions. Then
$$c_{\phi}=\dfrac{R^{'}_{\phi}(x)}{S_{\phi}(x)},$$ where
$R^{'}_\phi(x)$ is the differential of $R_\phi(x)$. Especially, for
any positive integer $n$, the map $[n]$ on $E$ is an endomorphism.
Set $[n](x,y)=(R_n(x), yS_n(x))$. From $c_{[1]}=1$ and
$[n]=[1]+[(n-1)]$, we have
$$c_{[n]}=\frac{R^{'}_n(x)}{S_n(x)}=n.$$

For any $Q=(x_Q,y_Q)\in E$, and any $$P=(x_P,y_P)\in
\Lambda=\{P=(x_P,y_P)\in E(\overline{K})\mid nP=Q\},$$ we have
$y_P=\frac{y_Q}{S_n(x_P)}$. Therefore, Theorem \ref{meanvalue} gives
$$\sum\limits_{P\in \Lambda}\frac{1}{S_n(x_P)}=\sum\limits_{P\in
\Lambda}\frac{y_P}{y_Q}=\frac{1}{y_Q}\sum\limits_{P\in
\Lambda}y_P=n^3.$$ Thus $$\sum\limits_{P\in
\Lambda}\frac{1}{R^{'}_n(x_P)}=\sum\limits_{P\in
\Lambda}\frac{1}{n\cdot S_n(x_P)}=\frac{1}{n}\sum\limits_{P\in
\Lambda}\frac{1}{S_n(x_P)}=n^2,$$ and
$$\sum\limits_{P\in
\Lambda}\frac{x_Q}{R^{'}_n(x_P)}=x_Q\sum\limits_{P\in
\Lambda}\frac{1}{R'_n(x_P)}=n^2x_Q=\sum\limits_{P\in \Lambda}x_P.$$

\end{document}